\input amssym.tex

\scrollmode
\nopagenumbers

\font\Bbb=msbm10

\def\1{\'{\i}} 

\def\N{\hbox{\Bbb N}} 
\def\R{\hbox{\Bbb R}} 
\def\qed{\quad\strut\vrule width.25cm height3.5pt \vskip5pt}

\centerline{\bf A notion of primitive that fits the Riemann integral}
\medskip

\rightline{Winston Alarc\'on-Athens}

\rightline{Escuela de Matem\'atica, U.C.R.}

\bigskip
{\narrower\narrower\bigskip\noindent{\bf Abstract}. We present a notion of primitive which corresponds exactly with the Riemann integral. We obtain a characterization of the integrability in the sense of Riemann which produces a Fundamental Theorem of Calculus without special assumptions. We prove a lemma which generalizes the theorem on functions with zero derivative in an interval. We apply these concepts, sketching a simplified alternative presentation of the Riemann integral in \R.\bigskip}

\noindent
{\bf 1. Introduction}. As it is well known, the integration in the sense of Riemann and the antiderivation are far from being equivalent operations and none of them reduces to the other: there are real functions differentiable at all points in a closed interval $ [a, b]  $ whose derivative is bounded but not Riemann  integrable on $ [a, b] $ ([{\bf 3}]).  Thus, for the fulfillment of the equality
$$
\int_a^b f(x)\,dx = F(b) - F(a)
$$
it is not sufficient to assume that $ F $ is an antiderivative of the bounded function $ f $  in $ [a, b] $: we must also require that $ f $ is Riemann--integrable on $ [a, b] $.

\smallskip
Moreover, there are functions that are integrable in the sense of Riemann on an interval $ [a, b] $ (as the characteristic function of the Cantor set in $ [0,1] $), whose integral from $ a $ to $ x $ defines a function of the variable $ x $ which is not a primitive from the first, at least in the usual sense.

\smallskip
Of course, this means that the Fundamental Theorem of Calculus requires special assumptions when working with the Riemann integral.

\smallskip
In what follows we will discuss an elementary generalization of the notion of  primitive or antiderivative which corresponds exactly with the Riemann integral, producing a TFC without special assumptions, along with a simplified presentation of the Riemann integral in \R.

\bigskip\noindent
{\bf  Pre-primitives and Riemann's primitives}. As we know that differenciation is not the inverse operation of Riemann integration, we renounce  to consider as a central issue the limit $$
\lim_{y \to x} {F(y) - F(x) \over y - x} \,; 
$$ 
instead, we will impose certain requirements on the ratio of increments that appears below the previous limit, thus establishing a first stage of the generalization we seek.

\smallskip
In what follows, $ J $ denote a non-degenerate interval of the real line \R\ and, for brevity, we say that a function $ f \colon J  \to \R $ is {\it bounded on compacts \/} if the restriction of $ f $ to any compact subinterval of $ J $ is bounded.

\bigskip\noindent
{\bf Definition A}. Let $ f \colon J \to \R $ be a bounded on compacts. We say that a function $ F \colon J \to \R $ is a {\it pre-primitive\/} of $ f $ in $ J $ if given $ x $, $ y \in J $ with $ x \neq y $, the ratio of increments
$$ F (y) - F (x) \over y - x $$
is between each pair of lower and upper bounds for $ f $ in the compact subinterval $ J $ of endpoints $ x $, $ y $.

\bigskip
First checks that Definition A is a generalization of the usual notion of primitive or antiderivative.

\bigskip\noindent
{\bf Theorem 1}.  Let $ f \colon J \to \R $ bounded on compacts. The primitives of $ f $ in $ J $ (if any) are pre-primitives of $ f $ in $ J $.

\medskip
{\bf Proof}. If $ F \colon J \to \R $ is a primitive of $ f $ in $ J $, let $ x $, $ y $ be two distinct points of $ J $, say $ x <y $. Since $ F $ is continuous on $ [x, y] $ and differentiable in $ \,] x, y [\, $, by the Mean Value Theorem for derivatives, there exists $ z \in \,] x, y [\, $ such that
$$
f(z)=F'(z) = {F(y) - F(x) \over y - x}
$$
Now, if $ p $, $ q $ are any lower and upper bounds for $ f $ in $ [x, y] $, we have $ p \leq f (z) \leq q $, thus:
$$
p \leq {F(y) - F(x) \over y - x} \leq q \,.
$$
This prove the primitive $ F $ is a pre-primitive of $ f $ in $ J $.\qed

\bigskip
Pre-primitives have some properties that resemble primitives. The next two theorems describe this situation.

\bigskip\noindent
{\bf Theorem 2}. Pre-primitives are continuous. In fact, they possess the Lipschitz property in each compact subinterval of its domain.

\medskip
{\bf Proof}. If $ F $ is a pre--primitive of a bounded on compacts function $ f $ in an interval $ J $ and if $ a $, $ b \in J $ satisfy $ a < b $, let $ L = \max \big \{| m |, | M | \big \} $, where $ m $ and $ M $ respectively denote the infimun and the supremun of $ f $ in $ [a, b]$. If $ x $, $ y $ are two distinct points of $ [a, b] $ with $ x < y $, we have:
$$
-L \leq m \leq {F(y) - F(x) \over y - x} \leq M \leq L \,,
$$
from which we get:
$$
\big| F(y) - F(x) \big| \leq L \, | y-x | \,.
$$
This shows $ F $ is lipschitzian in each compact subinterval $ J $, which in turn implies that $ F $ is continuous.\strut\qed

\bigskip\noindent
{\bf Theorem 3}. If $ F $ is a pre--primitive of $ f $ in the interval $ J $ and $ f $ has right lateral limit $ f(x^+) $ (respectively, left lateral limit $ f(x^-) $) at $ x \in J $, then $ F $ has  right derivative $ D^+ F(x) $ (respectively, left derivative $ D^-F(x) $) in $ x $, verifying their equality. In particular, if $ f $ is continuous at $ x $ then $ F $ is differentiable at $ x $, verifying $ F'(x) = f(x) $. If $ f $ is continuous, $ F $ is a primitive or antiderivative of $ f $ in $ J $.

\medskip
{\bf Proof}. Let $ y \in J $ with $ x < y $. Let $ m(y) $ and $ M(y) $ the infimum and the supremum of $ f $ in $ [x, y] $. By hypothesis we have:
$$
m(y)\leq{F(y)-F(x)\over y-x}\leq M(y)\,.
$$
Suppose that $ f $ has right lateral limit $ f (x^+) $ at $ x $, then:
$$
\lim_{y\to x^+}m(y) = \lim_{y\to x^+}M(y) = \lim_{y\to x^+}f(y) = f(x^+)\,.  
$$
This shows that
$$
D^+F(x)=\lim_{y\to x^+}{F(y)-F(x)\over y-x}=f(x^+)\,.
$$
Similarly it is proved that $ D^-F(x) = \lim_{y \to x^-} = f(x^-) $ when $ f $ has the left lateral limit $ f (x^{-} ) $ in $ x $.\qed

\bigskip
A nice aspect of the pre--primitives is that they certainly exist.

\bigskip\noindent
{\bf Theorem 4}. If $ f \colon J \to \R $ is bounded on compacts, then $ f $ has pre--primitives at $ J $.

\medskip
{\bf Proof}. 
Let $ c \in J $. Define the $ \underline{F}_c \colon J \to \R $ by:
$$
\underline{F}_c(x)=\cases{\underline{\cal I}(f;c,x)&si $x>c$\cr
\qquad 0&si $x=c$\cr
-\underline{\cal I}(f;x,c)&si $x<c$}
$$ 
where for all $ a $, $ b \in J $ with $ a <b $, $ \underline {\cal I}(f, a, b) $ denotes the supreme of the set $ \underline {\cal S}(f, a, b) $ of the {\it lower sums\/} for $ f $ in $ [a, b] $, ie sums of the form $ \sum_{i = 1}^n p_i (x_i - x_{i-1}) $, $ n \in \N $, $ a = x_0 \leq x_1 \leq \ldots \leq x_n = b $ and $ p_i $ a lower bound for $ f $ in $ [x_{i-1}, x_i] $, $ i = 1, \ldots, n $. Noting that if $ x <y <z $, we have:
$$\underline{\cal S}(f;x,z) = \big\{ s + s'\,\big|\, s\in \underline{\cal S}(f;x,y),\, s'\in \underline{\cal S}(f;y,z)\big\}$$
and we can use the additive property of the supreme [{\bf 1}, page~12], obtaining:
$$
\underline{\cal I}(f;x,y) + \underline{\cal I}(f;y,z) = \underline{\cal I}(f;x,z)\,.
$$
This allows to write:
$$
\underline{F}_c(y) - \underline{F}_c(x) = \underline{\cal I}(f;x,y)
$$
if $ x <y $. Now, if $ p $, $ q $ are any lower and upper bounds for $ f $ in $ [x, y] $, the definition of $ \underline {\cal I} (f, x, y) $ implies:
$$
p\,(y-x) \leq \underline{\cal I}(f;x,y) \leq q\,(y-x)\,,
$$
and combining these results we obtain:
$$
p \leq {\underline{F}_c(y) - \underline{F}_c(x) \over y-x}  \leq q\,.
$$
This proves $ \underline{F}_c \colon J \to \R $ is a pre--primitive of $ f $ in $ J $.\qed

\bigskip
An important difference between the notions of primitive and pre--primitive is that the difference of two pre--primitives of a function in an interval {\it is not necessarily constant}. To do this, simply consider the characteristic function of the rationals (or Dirichlet's function) $ \chi_{\Bbb Q} $. Since $  \chi_{\Bbb Q} $ is 1 in every rational number and is 0 in all other real numbers, then if $ x <y $ and $ p $, $ q \in \R $ are any lower and upper bounds of $ \chi_{\Bbb Q} $ in $ [x, y] $, we have: $ p \leq 0 <1 \leq q $. Consequently, any function $ F \colon \R \to \R $ in which the secants to its graph have non-negative slope no greater than 1, is a pre--primitive of $  \chi_{\Bbb Q} $ in \R. For example, the functions $ F_1(x) = 0 $ and $ F_2 (x) = x $ is a pair of such pre--primitives of $ \chi_{\Bbb Q} $. Since de difference $ F_1 - F_2 $ is not constant, then at least one of them is not  a primitive of $  \chi_{\Bbb Q}$ in \R, in any usual sense.

\medskip
However, there are important classes of functions whose pre--primitives differ in constants. For example, the pre--primitives of a continuous function on the interval $ J $ are, by Theorem 3, antiderivatives of $ f $ in $ J $, so they differ in constants. The same Theorem 3, together with the following lemma, will provide an even broader class of functions bounded on compacts whose pre--primitives differ in constants.

\bigskip\noindent
{\bf Lemma}. Let $ H \colon J \to \R $ such that at each point $ x\in J $, $ H $ has right  derivative (resp., left  derivative) identically null. Then $ H $ is constant.

\medskip
{\bf Proof}. 
Suppose $ H $ has null right  derivative at each point of the interval $ J $. Let us further assume, by contradiction, that $ H $ is not constant in $ J $. Then exists $a$, $b \in J $, say $a<b$, such that $ H(a) \neq H(b) $. Let $ L = | H (b)-H (a) |> 0 $. Consider the set
$$
{\cal C}= \bigg\{x\in\,]a,b]\,\bigg|\,{|H(x)-H(a)|\over x-a}\leq {L\over 2(b-a)}\bigg\}
$$
Since $ {L \over 2 (b-a)}> 0 $ and $ \lim_{x \to a^+} {H (x)-H (a) \over x-a} = 0 $, the set $ \cal C $ is not empty. Let $ c $ the supremum of $ \cal C $. Since $ \Big| {H (b)-H (a) \over b-a} \Big| = {L \over b-a}> {L \over2 (b-a)} $, we have $ a <c <b $. Let
$$
{\cal D}=\bigg\{x\in\,]c,b]\,\bigg|\, {|H(x)-H(c)|\over x-c}\leq{L\over 2(b-a)}\bigg\}\,.
$$
${\cal D}\neq\emptyset$, since $\lim_{x\to c^+}{H(x)-H(c)\over x-c}=0$. Let $ d \in \cal D $. We have: $ c <d $. 

\smallskip
\noindent
However
$$\eqalign{
\big|H(d)-H(a)\big|&=\big|H(d)-H(c)+H(c)-H(a)\big|\cr
&\leq \big|H(d)-H(c)\big| + \big|H(c)-H(a)\big|\cr
&={|H(d)-H(c)|\over d-c}(d-c)+{|H(c)-H(a)|\over c-a}(c-a)\cr
&\leq {L\over 2(b-a)}\Big((d-c)+(c-a)\Big)={L\over 2(b-a)}(d-a)\,.
}$$
As a result:
$$
{|H(d)-H(a)|\over d-a}\leq{L\over2(b-a)}
$$
and hence it follows that $ d \in \cal C $. We have obtained the contradiction $d\leq c<d$. Thus it is proved that $ H $ is constant $ J $. The case where $ H $ has null left derivative at each point of the interval $ J $ is proved similarly.\qed

\bigskip\noindent
{\bf Theorem 6}. Let $ f \colon J \to \R $ bounded on compacts. A sufficient condition for the pre--primitives of $ f $ in $ J $ to differ in constants, is that $ f $ has right  lateral  limit  (respectively, left lateral limit) at each point of $ J $.

\medskip
{\bf Proof}.  If $ f $ has right lateral limit at each point of $ J $, then by Theorem 3, if $ F $ and $ G $ are pre--primitives of $ f $ in $ J $, both have right derivative at every point of $ J $ and $ D^+ F(x) = D^+ G(x) = f(x^+)$. This shows  the function $ H = F-G $ has  right derivative, null at each point of $ J $ and Lemma 5 implies that $ H $ is constant $ J $.\qed

\bigskip
This result motivates the following definition:

\bigskip\noindent
{\bf Definition B}. Let $ f \colon J \to \R $ bounded on compacts. The pre--primitives of $ f $ in $ J $  will be called {\it Riemann's primitives\/} (or $ \cal R $--primitives) of $ f $ in $ J $  if they differ  in constants.

\bigskip
The importance of Definition B lies in the following theorem, whose part ``if''  is a TFC (part II) without special assumptions:

\bigskip\noindent
{\bf Theorem 7}. Let $ f \colon J \to \R $ bounded on compacts. Then $ f $ is Riemann integrable on the compact subintervals of $ J $ if and only if the pre--primitives of $ f $ in $ J $ are Riemanns's primitives of $ f $ in $ J $. If $ f $ is Riemann integrable on the compacts subintervals of $ J $ and $ a $, $ b \in J $, then:
$$
\int_a^b f(x)\,dx = F (b) - F (a)\,,
$$
where $ F $ is any $ \cal R $--primitive of $ f $ in $ J $.

\medskip
{\bf Proof}.  As we saw in the proof of Theorem 4, the function $ \underline {F}_c $ defined there is a pre--primitive of $ f $ in $ J $. Similarly it can be proved that the following function $\overline{F}_c\colon J\to\R$ is a pre--primitive of $ f $ in $ J $:
$$
\overline{F}_c(x)=\cases{\overline{\cal I}(f;c,x)&si $x>c$\cr
\qquad 0&si $x=c$\cr
-\overline{\cal I}(f;x,c)&si $x<c$}
$$ 
where for all $ a $, $ b \in J $ with $ a <b $, $ \overline{\cal I} (f; a, b) $ denotes the infimum of all the {\it superior sums\/} for $ f $ in $ [a, b] $, ie sums of the form $ \sum_{i = 1}^n q_i (x_i - x_{i-1}) $ with $ n \in \N $, $ a = x_0 \leq x_1 \leq \ldots \leq x_n = b $ and $ q_i $ any upper bound of $ f $ in $ [x_{i-1}, x_i] $, $ i = 1, \ldots, n $.

\smallskip
i) 
Suppose the pre--primitives of $ f $ in $ J $ are $ \cal R $--primitives of $ f $ in $ J $. Then the functions $ \underline{F}_c $ and $ \overline{F}_c $ differ by a constant and that constant is 0 since $ \underline{F}_c (c) = 0 = \overline{F}_c (c) $. Consequently $ \underline{F}_c = \overline{F}_c $. Now, if $ a $, $ b \in J $ are such that $ a <b $, we have:
$$
\underline{\cal I}(f;a,b)=\underline{F}_c(b)-\underline{F}_c(a)=\overline{F}_c(b)-\overline{F}_c(a)=\overline{\cal I}(f;a,b)\,,
$$
This proves that $ f $ is Riemann integrable on $ [a, b] $.

\smallskip
ii) 
Now suppose that $ f $ is Riemann integrable on the compact subintervals of $ J $ and $ a $, $ b \in J $ are such that $ a <b $. Let $ a = x_0 <x_1 <x_2 <\ldots <x_n = b $ an arbitrary partition of the interval $ [a, b] $. For each $ i = 1, \ldots, n \,$ let $ p_i $ and $ q_i $ be any lower and upper bounds for $ f $ in the interval $ [x_{i-1}, x_i] $. If $ F $ is any pre--primitive of $ f $ in $ J $, for each $ i = 1, \ldots, n $ occurs:
$$
p_i \leq {F(x_i)-F (x_{i-1}) \over x_i - x_{i-1}} \leq q_i \,.
$$
Multiplying all by $ x_i-x_{i-1}> 0 $ and adding member to member from $ i = $ 1 to $ i = n $, we obtain:
$$
\sum_{i = 1}^n p_i \, (x_i - x_{i-1}) \leq F (b) - F (a) \leq \sum_{i = 1}^n q_i \, (x_i - x_{i-1})
$$
because the telescopic sum $ \sum_{i = 1}^n F(x_i)-F(x_{i-1}) $ collapses to the expression $F (b) - F (a) $.

\smallskip
We have proved that if $ F $ is any pre--primitive of $ f $ in $ J $ and $ a $, $ b \in J $ are such that $ a <b $, then the difference $ F (b)-F (a) $ is at the same time, an upper bound of the lower sums for $ f $ in $ [a, b] $, and a lower bound of the upper sums for $ f $ in $ [a, b] $. Therefore, by definition of $ \underline{\cal I}(f, a, b) $ and $ \overline{\cal I} (f, a, b) $, we have:
$$
\underline{\cal I}(f;a,b) \leq F(b) - F(a) \leq \overline{\cal I}(f;a,b)\,.
$$
As we are assuming that $ f $ is Riemann integrable on the compact subintervals of $ J $, we conclude that:
$$
\int_a^b f(x)\,dx = F(b) - F(a)\,,
$$
for any  pre--primitive $ F $ of $ f $ in $ J $, where now we can lift the restriction $ a <b $ and require only that $ a $, $ b \in J $. This allows to conclude that if $ F $ and $ G $ are any two pre--primitives of $ f $ in $ J $, then $ F (x)-F (a) = G (x)-G (a) $ for all $ x \in J $, which shows that the difference $ F-G $ is constant in $ J $. This proved that the pre--primitives of $ f $ in $ J $ are Riemann's primitives of $ f $ in $ J $, in the sense of Definition B.\qed

\bigskip\noindent
{\bf Theorem 8}. Let $ f \colon J \to \R $ bounded on compacts. A sufficient condition for $ f $ to be Riemann integrable on the compact subintervals of $ J $ is that $ f $ has  limit from the right at each point of $ J $, or  limit from the left at each point of $ J $. [{\bf 2}].

\medskip
{\bf Proof}.  An immediate consequence of Theorems 6 and 7.\qed

\bigskip
As direct consequences of the theorems just discussed, we obtain the following classical results:

\bigskip\noindent
{\bf Corollary 9}. If $ f \colon J \to \R $ is monotone, then $ f $ is Riemann integrable on the compact subintervals of $ J $.

\medskip
{\bf Proof}.  If $ f \colon J \to \R $ is monotone, $ f $ is bounded on compacts and at every point of $ J $, $ f $ has limit from the left  if $ f $ is increasing or has limit from the right if $ f $ is decreasing. In either cases, by Theorem 8 $ f $ is Riemann--integrable on the compact subintervals of $ J $.\qed

\bigskip\noindent
{\bf Corollary 10}. Any continuous real function on an interval $ J $ admits primitives in $ J $ and is Riemann integrable on the compact subintervals of $ J $.

\medskip
{\bf Proof}.  If $J\subseteq\R$ is an interval y $ f \colon J \to \R $ is continuous, $ f $ is bounded on compacts and satisfy the hypotheses of Theorem 8.\qed

\medskip
Corollary 10 is also a direct consequence of Theorems 3, 4 and 7, without the intervention of Theorem 6 (nor of Lemma 5 that is its base). In fact: By Theorem 4, $ f $ has pre--primitives in $ J $ and Theorem 3 such pre--primitives are antiderivatives of $ f $ in $ J $, so the pre--primitives of $ f $ in $ J $ are $ \cal R $--primitives of $ f $ in $ J $. Then, by Theorem 7, $ f $ is Riemann integrable on compact subintervals of $ J $.\qed

\medskip
Note that in the two previous proofs of Corollary 10, the property of {\it uniform continuity\/} of $ f $ on the compacts subintervals of $J$ has not been used.

\bigskip\noindent
{\bf Corollary 11. (Fundamental Theorem of Calculus, Part I)}. Let $ f \colon J \to \R $ Riemann integrable on compact subintervals of $ J $. Let $ c \in J $ and consider the function $ F_c \colon J \to \R $ defined by:
$$
F_c (x) = \int_c^x f (t) \, dt \,.
$$
Then $ F_c $ is continuous and if $ f $ is continuous at $ x \in J $, $F_c $ is differentiable at $ x $, verifying $ F '(x) = f (x) $.

\medskip
{\bf Proof}. By Theorem 7, $ F_c (x)  = G(x)-G(c) $, where $ G $ is any pre--primitive of $ f $ in $ J $. Since $ F_c $ and $ G $ differ by a constant,  $ F_c$ is a  pre--primitive of $ f $ in $ J $. By Theorem~3, $ F_c $ is continuous on $ J $, and if $ f $ is continuous at $ x \in J $, $F_c $ is differentiable in $ x $, verifying $ F'_c (x) = f (x)$.\qed

\bigskip\noindent
{\bf Corollary 12}. (Fundamental Theorem of Calculus, Part II). Let $ a $, $ b \in \R $ with $ a <b $. Let $ f \colon [a, b] \to \R $, a bounded function. If $ f $ is Riemann integrable on $ [a, b] $ and if $ f $ has primitives or antiderivatives in $ [a, b] $, then
$$
\int_a^b f(x)\,dx = F (b) - F (a) \,,\eqno(7)
$$
where $ F $ is any primitive of $ f $ on $ [a, b] $.

\medskip
{\bf Proof}. Let $ F $ an arbitrary primitive $ f $ on $ [a, b] $. By Theorem~1, $ F $ is a pre--primitive of $ f $ in $ [a, b] $. On the other hand, since $ f $ is Riemann integrable on $ [a, b] $, then $ f $ is Riemann integrable on each closed subinterval of $ [a, b] $ and Theorem 7 assert the equality.\qed

\bigskip\noindent
{\bf References.}

\medskip\noindent
{\bf[1]} T.~M.~Apostol: {\it An\'alisis Matem\'atico}. Revert\'e, 2a ed., Barcelona, {\oldstyle 1977}.

\smallskip\noindent
{\bf[2]} R.~C.~Metzler: {\it On Riemann integrability}. AMM vol.~78, n.~10, p.~1129.

\smallskip\noindent
{\bf[3]} I.~Muntean: {\it Clasificaci\'on de algunas funciones reales en un intervalo compacto}. Ciencias Matem\'aticas, vol.~1, n.~1, p\'ag.~39, Universidad de Costa Rica, San Jos\'e, {\oldstyle 1990}.

\bye